\documentclass[12pt]{article}      %{article} was 12pt latex e
\usepackage{amssymb}
\usepackage{eucal,amsthm}
\newtheorem{theorem}{Theorem}[section]
\newtheorem{prop}[theorem]{Proposition}

\newtheorem{lemma}[theorem]{Lemma}
\newtheorem{coro}[theorem]{Corollary}
\newtheorem{prop-def}{Proposition-Definition}[section]

\newtheorem{exam}[theorem]{Example}
\newcommand{\nc}{\newcommand}

\newcommand{\delete}[1]{}
\newcommand{\mlabel}[1]{\label{#1}}  % Use this line to suppress names
\nc{\dfootnote}[1]{} %{\footnote{#1}}

\nc{\bin}[2]{ (_{\stackrel{\scs{#1}}{\scs{#2}}})}  %binomial coeff
\nc{\binc}[2]{ \left (\!\! \begin{array}{c} \scs{#1}\\
    \scs{#2} \end{array}\!\! \right )}  %binomial coeff 
\nc{\bincc}[2]{  \left ( {\scs{#1} \atop
    \vspace{-1cm}\scs{#2}} \right )}  %binomial coeff  
\nc{\bs}{\bar{S}}
\nc{\Cq}{{C\!\!<\!\! {\mathfrak q}\!\!>}}
\nc{\la}{\longrightarrow}
\nc{\rar}{\rightarrow}
\nc{\dar}{\downarrow}
\nc{\dap}[1]{\downarrow \rlap{$\scriptstyle{#1}$}}
\nc{\uap}[1]{\uparrow \rlap{$\scriptstyle{#1}$}}
\nc{\defeq}{\stackrel{\rm def}{=}}
\nc{\dis}[1]{\displaystyle{#1}}
\nc{\dotcup}{\ \displaystyle{\bigcup^\bullet}\ }
\nc{\hcm}{\ \hat{,}\ }
\nc{\hts}{\hat{\otimes}}
\nc{\hcirc}{\hat{\circ}}
\nc{\lleft}{[}
\nc{\lright}{]}
\nc{\curlyl}{\left \{ \begin{array}{c} {} \\ {} \end{array}
    \right .  \!\!\!\!\!\!\!} 
\nc{\curlyr}{ \!\!\!\!\!\!\!
    \left . \begin{array}{c} {} \\ {} \end{array}
    \right \} }
\nc{\longmid}{\left | \begin{array}{c} {} \\ {} \end{array}
    \right . \!\!\!\!\!\!\!}
\nc{\ora}[1]{\stackrel{#1}{\rar}}
\nc{\ola}[1]{\stackrel{#1}{\la}}%${\Bbb Z}$
\nc{\scs}[1]{\scriptstyle{#1}}
\nc{\mrm}[1]{{\rm #1}}
\nc{\margin}[1]{\marginpar{}}   %{\rm #1}}
\nc{\dirlim}{\displaystyle{\lim_{\longrightarrow}}\,}
\nc{\invlim}{\displaystyle{\lim_{\longleftarrow}}\,}
\nc{\mvp}{\vspace{0.3cm}}
\nc{\tk}{^{(k)}}
\nc{\tp}{^\prime}
\nc{\ttp}{^{\prime\prime}}
\nc{\svp}{\vspace{2cm}}
\nc{\vp}{\vspace{8cm}}
\nc{\proofbegin}{\begin{proof}} % AMS command
\nc{\proofend}{\end{proof}} %AMS command
\nc{\modg}[1]{\!<\!\!{#1}\!\!>}
\nc{\intg}[1]{F_C(#1)}
\nc{\lmodg}{\!<\!\!}
\nc{\rmodg}{\!\!>\!}
\nc{\cpi}{\widehat{\Pi}}
\nc{\sha}{{\mbox{\cyr X}}}  %used to be \cyr
\nc{\shpr}{\diamond}    %Shuffle product
\nc{\vep}{\varepsilon}
\nc{\labs}{\mid\!}
\nc{\rabs}{\!\mid}

\nc{\ann}{\mrm{ann}}
\nc{\Aut}{\mrm{Aut}}
\nc{\can}{\mrm{can}}
\nc{\colim}{\mrm{colim}}
\nc{\Cont}{\mrm{Cont}}
\nc{\rchar}{\mrm{char}}
\nc{\cok}{\mrm{coker}}
\nc{\dtf}{{R-{\rm tf}}}
\nc{\dtor}{{R-{\rm tor}}}

\nc{\Div}{{\mrm Div}}
\nc{\End}{\mrm{End}}
\nc{\Ext}{\mrm{Ext}}
\nc{\Fil}{\mrm{Fil}}
\nc{\Frob}{\mrm{Frob}}
\nc{\Gal}{\mrm{Gal}}
\nc{\GL}{\mrm{GL}}
\nc{\Hom}{\mrm{Hom}}
\nc{\hsr}{\mrm{H}}
\nc{\hpol}{\mrm{HP}}
\nc{\id}{\mrm{id}}
\nc{\im}{\mrm{im}}
\nc{\incl}{\mrm{incl}}
\nc{\length}{\mrm{length}}
\nc{\mchar}{\rm char}
\nc{\mpart}{\mrm{part}}
\nc{\ord}{\mrm{ord}}
\nc{\ql}{{\QQ_\ell}}
\nc{\qp}{{\QQ_p}}
\nc{\rank}{\mrm{rank}}
\nc{\rcot}{\mrm{cot}}
\nc{\rdef}{\mrm{def}}
\nc{\rdiv}{{\rm div}}
\nc{\rtf}{{\rm tf}}
\nc{\rtor}{{\rm tor}}
\nc{\res}{\mrm{res}}
\nc{\SL}{\mrm{SL}}
\nc{\Spec}{\mrm{Spec}}
\nc{\tor}{\mrm{tor}}
\nc{\Tr}{\mrm{Tr}}
\nc{\tr}{\mrm{tr}}

\nc{\ab}{\mathbf{Ab}}
\nc{\Alg}{\mathbf{Alg}}
\nc{\Bax}{\mathbf{Bax}}
\nc{\bfk}{{\bf k}}
\nc{\bfone}{{\bf 1}}
\nc{\base}[1]{{a_{#1}}}
\nc{\detail}{\marginpar{\bf More detail}
    \noindent{\bf Need more detail!}
    \svp}
\nc{\Diff}{\mathbf{Diff}}   
\nc{\gap}{\marginpar{\bf Incomplete}\noindent{\bf Incomplete!!}
    \svp}
\nc{\FMod}{\mathbf{FMod}}
\nc{\Int}{\mathbf{Int}}
\nc{\Mon}{\mathbf{Mon}}
\nc{\remarks}{\noindent{\bf Remarks: }}
\nc{\Rep}{\mathbf{Rep}}
\nc{\Rings}{\mathbf{Rings}}
\nc{\Sets}{\mathbf{Sets}}
\nc{\bill}[1]{\marginpar{\bf To Bill}\noindent{\bf To Bill:}
    {\tt #1}\\ }
\nc{\li}[1]{\marginpar{\bf To Li}\noindent{\bf To Li:}
    {\tt #1}\\ }

\nc{\BA}{{\Bbb A}}
\nc{\CC}{{\Bbb C}}
\nc{\DD}{{\Bbb D}}
\nc{\EE}{{\Bbb E}}
\nc{\FF}{{\Bbb F}}
\nc{\GG}{{\Bbb G}}
\nc{\HH}{{\Bbb H}}
\nc{\LL}{{\Bbb L}}
\nc{\NN}{{\Bbb N}}
\nc{\QQ}{{\Bbb Q}}
\nc{\RR}{{\Bbb R}}
\nc{\TT}{{\Bbb T}}
\nc{\VV}{{\Bbb V}}
\nc{\ZZ}{{\Bbb Z}}

\nc{\cala}{{\mathcal{ A}}}
\nc{\calc}{{\cal C}}
\nc{\cald}{\mathcal{D}}
\nc{\cale}{{\cal E}}
\nc{\calf}{{\cal F}}
\nc{\calg}{{\cal G}}
\nc{\calh}{{\cal H}}
\nc{\cali}{{\cal I}}
\nc{\call}{{\cal L}}
\nc{\calm}{{\cal M}}
\nc{\caln}{{\cal N}}
\nc{\calo}{{\cal O}}
\nc{\calp}{{\cal P}}
\nc{\calr}{{\cal R}}
\nc{\calt}{{\cal T}}
\nc{\calw}{{\cal W}}
\nc{\calx}{{\cal X}}
\nc{\CA}{\mathcal{A}}

\nc{\fraka}{{\frak a}}
\nc{\frakB}{{\frak B}}
\nc{\frakm}{{\frak m}}
\nc{\frakp}{{\frak p}}

\font\cyr=wncyr10

\vspace{-6cm}

\title{ 
Baxter Algebras and Hopf Algebras
}

\author{
George E. Andrews\\
Department of Mathematics\\ 
Pennsylvania State University\\ 
University Park, PA 16802, USA\\
(andrews@math.psu.edu)
\smallskip
\\
Li Guo\\
Department of Mathematics and Computer Science\\
Rutgers University at Newark\\ 
Newark, NJ 07102, USA\\
(liguo@newark.rutgers.edu)
\smallskip
\\
William Keigher\\
Department of Mathematics and Computer Science\\ 
Rutgers University at Newark\\
Newark, NJ 07102, USA\\
(keigher@newark.rutgers.edu)
\smallskip
\\
Ken Ono\\
Department of Mathematics\\ 
University of Wisconsin\\ 
Madison, WI 53706, USA\\
(ono@math.wisc.edu)
\thanks{
The first and fourth authors are supported by grants from the National
Science Foundation, and the fourth author is supported by Alfred P. Sloan,
David and Lucile Packard, and H. I. Romnes Fellowships.
}
}

\date{}

\begin{document}
\maketitle

\newpage

\begin{abstract}
By applying a recent construction of free Baxter algebras, 
 we obtain a new class of Hopf algebras
that generalizes the classical divided 
power Hopf algebra. 
We also study conditions under which these Hopf algebras are 
isomorphic. 
\end{abstract}

\section{Introduction}
\mlabel{intro}
Hopf algebras have their origin in Hopf's
seminal works on topological groups in the 1940s, and have
become fundamental objects
in many areas of mathematics and physics. 
For example, they are crucial 
to the study of algebraic groups, Lie groups, Lie algebras, and
quantum groups.  In turn, these areas have provided 
many of the  most important examples of Hopf algebras.

In this paper we construct new examples of Hopf algebras.
Our examples arise naturally in a combinatorial
study of Baxter algebras.
A Baxter algebra~\cite{Ba} 
is an algebra $A$ with a linear operator $P$ on $A$ 
that satisfies the identity 
\[ P(x)P(y) = P(xP(y)) + P(yP(x)) + \lambda P(xy) \]
for all $x$ and $y$ in $A$, where
$\lambda$, the  {\it weight}, is a fixed element in 
the ground ring of the algebra $A$.
Rota~\cite{Ro1} began a systematic study of Baxter 
algebras from an algebraic and combinatorial perspective
and suggested that they are related to
hypergeometric 
functions, incidence algebras 
and symmetric functions~\cite{Ro2,Ro3}. 
A survey of Baxter algebras with examples and applications can be found 
in \cite{Ro2,Ro3}, as well as in \cite{Gu1}. 

Free Baxter algebras were first constructed by 
Rota~\cite{Ro1} and Cartier~\cite{Ca} in the 
category of Baxter algebras with no identity 
(with some restrictions on the weight and the base ring).  
Recently, two of the authors~\cite{G-K1,G-K2} 
have constructed free Baxter algebras in a more 
general context including these classical constructions.
Their construction is in terms of {\it mixable shuffle 
products} which generalize the well-known 
{\it shuffle products} of path integrals 
as developed by Chen~\cite{Ch} and Ree~\cite{Re}. 
Here we show that a special case of 
the construction of these new Baxter algebras 
provides a large supply of new Hopf algebras.

The divided power Hopf algebra is one of the classical
examples of a Hopf algebra, and
it is not difficult to see that this algebra is the free Baxter algebra 
of weight zero on the empty set. 
The new Hopf algebras presented here generalize
this classical example.
In particular, we show that the
free Baxter algebra of
arbitrary weight on the empty set is a Hopf algebra. 

Here we describe the construction. 
%\dfootnote{This part can be moved to Section~\ref{sec:alg}.}
Let $C$ be a commutative algebra with identity $\bfone$, and
let $\lambda \in C$.
Define the sextuple $\cala=\cala_\lambda\defeq 
(A,\mu,\eta,\Delta,\vep,S)$, where 
%\newpage
\begin{eqnarray*}
&(i)\ \ \ \  A=A_\lambda =
    &\bigoplus_{n=0}^\infty C \base{n}, {\rm\ is\ the\ free\ }
    C\mbox{-module on the set\ } \{\base{n}\}_{n\geq 0}, \\
\\
&(ii)  \ \ \ \  \mu=\mu_\lambda:& A\otimes_C A \to A, \\
&&\displaystyle{\base{m}\otimes \base{n}\mapsto 
\sum_{k=0}^m \lambda^k \bincc{m+n-k}{m}\bincc{m}{k} 
    \base{m+n-k},} \\
&(iii)\ \ \ \ \  \eta=\eta_\lambda: & C\to A,\ 1\mapsto \base{0},\\
&(iv) \ \ \  \Delta=\Delta_\lambda:& A\to A\otimes_C A,  \base{n}\mapsto  
\sum_{k=0}^{n} \sum_{i=0}^{n-k} (-\lambda)^k \base{i}\otimes 
    \base{n-k-i},\\
&(v) \ \ \ \ \ \vep=\vep_\lambda: & A\to C,\ \base{n} \mapsto \left \{ 
\begin{array}{ll}
\bfone, & n=0,\\
\lambda\bfone, & n=1,\\
0, & n\geq 2
\end{array} \right . \\
&(vi)\ \ \ \ \  S=S_\lambda: &A\to A,\ \base{n} \mapsto 
  (-1)^n \sum_{v=0}^n \bincc{n-3}{v-3} \lambda^{n-v} \base{v}.
\mlabel{S}
\end{eqnarray*}
Here for any positive or negative integer $x$, $\bincc{x}{k}$ is defined 
by the generating function 
$\displaystyle{(1+z)^x=\sum_{k=0}^\infty \bincc{x}{k} z^k.}$
It was shown in~\cite{G-K1} that $(A_\lambda,\mu_\lambda,\eta_\lambda)$ 
is a Baxter algebra of weight $\lambda$ 
with respect to the 
operator
\[ P: A_\lambda\to A_\lambda,\ 
    \base{n} \mapsto \base{n+1}. \]
The main theorem in Section~\ref{sec:alg} is 

\begin{theorem}
For any $\lambda \in C$, 
$\cala_\lambda$ is a Hopf $C$-algebra. 
\mlabel{thm:alg}
\end{theorem}

When $\lambda=0$, we have the divided power Hopf algebra. In 
general, $\cala_\lambda$ will be called the {\it $\lambda$-divided power 
Hopf algebra}. 
The divided power algebra plays an important role in several areas of 
mathematics, including crystalline cohomology in number theory~\cite{B-O}, 
umbral calculus in combinatorics~\cite{R-R} and 
Hurwitz series in differential algebra~\cite{Ke}. 
We expect that the $\lambda$-divided power algebra $\cala_\lambda$ 
introduced here will play similar roles in 
these areas. 
To describe the application to umbral calculus, we recall that 
a sequence $\{p_n (x)\mid n\in\NN \}$ of polynomials in 
$C[x]$ is called a {\it sequence of binomial type} if 
\[ p_n(x+y)= \sum_{k=0}^n \bincc{n}{k} p_{n-k}(x) p_k(y)
\]
in $C[x,y]$.
The classic umbral calculus studies these sequences and their 
generalizations, such as Sheffer sequences and cross sequences. 
It is well-known that the theory of classic umbral calculus can be most conceptually 
described in the framework of Hopf algebras. 
Furthermore, most of the main results in umbral calculus follows from properties 
of the divided power Hopf algebra structure on the linear dual 
$\Hom_C(C[x],C)$, usually called the {\it umbral algebra} (see \cite{N-S,Rom} for details). 
This algebra is the completion of the divided power algebra equipped with an umbral shift 
operator which turns out to be a Baxter operator, making the umbral algebra 
the complete free Baxter algebra of weight zero on the empty set. See~\cite{Gu} 
for details. 
A program carried out there gives a generalization of the umbral calculus 
using the $\lambda$-divided power algebra. 

\smallskip 

It is natural to ask whether these Hopf algebras are isomorphic for 
different values of $\lambda$. We address this question in Section~\ref{sec:iso}
and obtain the following result. 

\begin{theorem}
Let $\lambda,\, \nu$ be in $C$. 
\begin{enumerate}
\item
If $(\lambda)=(\nu)$, then $\cala_\lambda$ and $\cala_\nu$ are 
isomorphic Hopf algebras. 
\mlabel{iso1}
\item
Suppose $C$ is a $\QQ$-algebra. If $\cala_\lambda$ and $\cala_\nu$ are 
isomorphic Hopf algebras, then $(\lambda)=(\nu)$. 
\mlabel{iso2}
\item
If $\cala_\lambda$ and 
$\cala_\nu$ are isomorphic Hopf algebras, then 
$\sqrt{(\lambda)}= \sqrt{(\nu)}$.
\mlabel{iso3}
\end{enumerate}
\mlabel{thm:iso}
\end{theorem}
We also prove a stronger but more technical version of the third statement 
in Proposition~\ref{pp:iso3}.  

\smallskip

According to Theorem~\ref{thm:iso}, for any $C$, there are at least two 
non-isomorphic $\lambda$-divided power Hopf algebras, namely $\cala_0$ and 
$\cala_1$. Furthermore, if $C$ is a $\QQ$-algebra, then there is a one-one 
correspondence between isomorphic classes of $\lambda$-divided power Hopf algebras 
and principal ideals of $C$. For other examples, see Example~\ref{ex:iso}.

\bigskip
\noindent
{\bf Acknowledgements} We thank Jacob Sturm for helpful discussions.

\section{On $\lambda$-divided power Hopf algebras}
\mlabel{sec:alg}

Since the constant $\lambda\in C$ will be fixed throughout this section,
we will drop the dependence on the subscript $\lambda$.
In section~\ref{ss:back} we recall the defining properties of a Hopf algebra.
The remainder of the section establishes that $\cala_\lambda$
satisfies these properties.
Since two of the authors proved that  
$(A,\mu,\eta)$ is a $C$-algebra~\cite{G-K1}, the proof
of Theorem 1.1 reduces to a step by step verification
of the defining properties of a Hopf algebra.

\subsection{Preliminaries}
\mlabel{ss:back}

Here we recall some basic definitions and facts for later 
reference.
All tensor products in this paper are taken over 
the fixed commutative ring $C$. 
Recall that a {\it cocommutative 
$C$-coalgebra} is a triple $(A,\Delta,\varepsilon)$ 
where $A$ is a $C$-module, $\Delta: A\to A\otimes A$ 
and $\varepsilon: A\to C$ are $C$-linear maps that make  
the following diagrams commute. 

\margin{coass}
\begin{equation}
%({\rm Coass})\ \ \ 
\begin{array}{ccc} 
    A & \ola{\Delta} & A\otimes A \\
    \dap{\Delta} && \dap{\id\otimes \Delta}\\
    A\otimes A &\ola{\Delta\otimes \id} & A\otimes A\otimes A
\end{array}
\mlabel{coass}
\end{equation}
\margin{coun}
\begin{equation}
%({\rm Coun}) \ \ \ 
\begin{array}{ccccc}
C\otimes A& \stackrel{\varepsilon\otimes \id}{\longleftarrow} & 
A\otimes A & \ola{\id\otimes \vep} & A\otimes C \\
&  {}_{\cong}\!\!\nwarrow &\uap{\Delta} & \nearrow_{\cong} \\
&& A&&
\end{array}
\mlabel{coun}
\end{equation}

\margin{cocom}
\begin{equation}
%({\rm Cocomm})\ \ \ 
\begin{array}{ccccc}
& & A & & \\
& {}^{\Delta} \!\!\swarrow & & \searrow^{\Delta}& \\
A\otimes A & & \ola{\tau_{A,A}} & &A\otimes A
\end{array}
\mlabel{cocom}
\end{equation}
where $\tau_{A,A}: A\otimes A\to A\otimes A$ is defined 
by $\tau_{A,A}(x\otimes y)=y\otimes x$. 
The $C$-algebra $C$ has a natural structure of 
a $C$-coalgebra 
with 
\[ \Delta_C: C\to C\otimes C,\ c\mapsto c\otimes 1,\ c\in C\]
and 
\[ \vep_C=\id_C: C\to C. \] 
We also denote the multiplication in $C$ by 
$\mu_C$. 

Recall that a $C$-{\it bialgebra} is a quintuple 
$(A,\mu,\eta,\Delta,\vep)$ where $(A,\mu,\eta)$ 
is a $C$-algebra and $(A,\Delta,\vep)$ is a $C$-coalgebra 
such that $\mu$ and $\eta$ are morphisms of $C$-coalgebras. 
In other words, we have the commutativity of the following 
diagrams. 

\margin{md}
\begin{equation}
\begin{array}{ccc}
 A\otimes A & \ola{\mu} & A \\
\llap{$\scriptstyle{(\id\otimes\tau\otimes \id)(\Delta\otimes 
\Delta)}$}
\downarrow
&&\dap{\Delta}\\
A\otimes A\otimes A\otimes A &\ola{\mu\otimes \mu} &
    A\otimes A
\end{array}
\mlabel{md}
\end{equation}

\margin{me}
\begin{equation}
\begin{array}{ccc}
A\otimes A & \ola{\vep\otimes \vep} & C\otimes C\\
\dap{\mu} && \dap{\mu_C}\\
A & \ola{\vep} & C
\end{array}
\mlabel{me}
\end{equation}

\margin{ed}
\begin{equation}
\begin{array}{ccc}
C &\ola{\eta} & A \\
\dap{\Delta_C} && \dap{\Delta}\\
C\otimes C &\ola{\eta\otimes \eta} & A\otimes A
\end{array}
\mlabel{ed}
\end{equation}

\margin{ee}
\begin{equation}
\begin{array}{ccccc}
C && \ola{\eta} && A\\
& {}_{\id}\searrow && \swarrow_{\vep} & \\
&& C &&
\end{array}
\mlabel{ee}
\end{equation}

Let $(A,\mu,\eta,\Delta,\vep)$ be a $C$-bialgebra. 
For $C$-linear maps $f,\, g:A\to A$, the convolution 
$f\star g$ of $f$ and $g$ is the composition of the maps 
\[ A\ola{\Delta} A\otimes A \ola{f\otimes g} A\otimes A 
    \ola{\mu} A.\]
A $C$-linear endomorphism $S$ of $A$ is called an {\it antipode} 
for $A$ if 
\margin{anti}
\begin{equation}
 S\star \id_A = \id_A\star S =\eta\circ \vep. 
\mlabel{anti}
\end{equation}
A {\it  Hopf algebra} is a bialgebra $A$ with an antipode $S$.

\subsection{Coalgebra Properties}
\mlabel{coal}

We verify that $(A,\Delta,\vep)$ satisfies the axioms of 
a coalgebra characterized by the diagrams~(\ref{coass}), 
(\ref{coun}) and (\ref{cocom}).

%\subsection{Commutativity of diagram (\ref{coass})}
To prove (\ref{coass}), we only need to verify that, 
for each $n\geq 0$, 
\margin{ca}
\begin{equation}
(\Delta\otimes \id)(\Delta(\base{n})) = 
(\id\otimes \Delta)(\Delta(\base{n})). 
\mlabel{ca}
\end{equation}
Unwinding definitions on the left hand side we get 
$$ (\Delta\otimes \id)(\Delta(\base{n})) 
= \sum_{k=0}^{n}\sum_{i=0}^{n-k}\sum_{\ell=0}^i 
\sum_{j=0}^{i-\ell} (-\lambda)^{k+\ell} \base{j}\otimes 
\base{i-\ell-j}\otimes \base{n-k-i}.$$
Exchanging the third and the fourth summations, and then 
exchanging the second and the third summations, we obtain  
$$ (\Delta\otimes \id)(\Delta(\base{n})) 
=\sum_{k=0}^n\sum_{j=0}^{n-k}\sum_{i=j}^{n-k}\sum_{\ell=0}^{i-j} 
(-\lambda)^{k+\ell} \base{j}\otimes \base{i-\ell-j}\otimes 
    \base{n-k-i}.$$
Replacing $i$ by $i+j$ gives us 
$$(\Delta\otimes \id)(\Delta(\base{n})) 
=\sum_{k=0}^n\sum_{j=0}^{n-k}\sum_{i=0}^{n-k-j}\sum_{\ell=0}^i 
(-\lambda)^{k+\ell} \base{j}\otimes \base{i-\ell}\otimes 
    \base{n-k-i-j}.$$
Exchanging the third and the fourth summations, followed by 
a substitution of $i$ by $i+\ell$, gives us 
$$ (\Delta\otimes \id)(\Delta(\base{n})) 
= \sum_{k=0}^n \sum_{j=0}^{n-k}\sum_{\ell=0}^{n-k-j} 
\sum_{i=0}^{n-k-j-\ell} (-\lambda)^{k+\ell} \base{j}\otimes 
\base{i}\otimes \base{n-k-i-\ell-j}.$$
We get the same expression after unwinding 
definitions on the right hand side. 
This proves equation~(\ref{ca}). 

%\subsection{Commutativity of diagram~(\ref{coun})}
Before proving the commutativity of diagram~(\ref{coun}) 
we display a lemma. 

\begin{lemma}
For any integers $n$ and $\ell$ with 
$n\geq \ell\geq 0$, we have 

\margin{cu2}
\begin{equation}
\sum_{k=0}^{n-\ell} (-\lambda)^k \vep (\base{n-\ell-k}) 
    =\left \{ \begin{array}{ll} \bfone, & \ell=n, \\
    0, & \ell< n. \end{array} \right . 
\mlabel{cu2}
\end{equation}
\mlabel{lem:cu}
\end{lemma}
\proofbegin
When $\ell=n$, we have 
$$ \sum_{k=0}^{n-\ell} (-\lambda)^k \vep (\base{n-\ell-k}) 
= \vep (\base{0}) = \bfone. $$
When $\ell< n$, we have 
$$
\sum_{k=0}^{n-\ell} (-\lambda)^k \vep (\base{n-\ell-k}) 
= (-\lambda)^{n-\ell-1} \vep(\base{1}) + 
    (-\lambda)^{n-\ell} \vep(\base{0}).  $$
By the definition of $\vep(\base{n})$, the right hand side is 
$$ (-\lambda)^{n-\ell-1} \lambda \bfone+ 
(-\lambda)^{n-\ell} \bfone= 0.$$
\proofend

We can now prove (\ref{coun}). 
Consider the left triangle in (\ref{coun}). 
For each $n\geq 0$, we have 

$$(\vep\otimes \id)(\Delta (\base{n})) 
= \sum_{k=0}^n \sum_{i=0}^{n-k} (-\lambda)^k \vep(\base{i}) 
\otimes \base{n-k-i}.$$
A substitution $i=n-k-\ell$ and then an exchange of the order 
of summations give us
$$(\vep\otimes \id)(\Delta (\base{n})) 
= \sum_{\ell=0}^n (\sum_{k=0}^{n-\ell} (-\lambda)^k 
    \vep (\base{n-k-\ell}))\otimes \base{\ell}$$
which is $\bfone\otimes \base{n}$ by Lemma~\ref{lem:cu}. 
This proves the commutativity of the left triangle 
in (\ref{coun}). 

The proof of the right triangle in (\ref{coun}) is similar. 

The cocommutativity of diagram~(\ref{cocom}) is easy to verify: 

$$\tau_{A,A}(\Delta(\base{n})) 
= \sum_{k=0}^n \sum_{i=0}^{n-k} (-\lambda)^k 
    \base{n-k-i}\otimes \base{i}.$$ 
After replacing $i$ by $n-k-j$, we see that it is the same as 
$\Delta(\base{n})$. 

Hence we have shown that $(A,\Delta,\vep)$ is a 
cocommutative $C$-coalgebra.

\subsection{Compatibility}
\mlabel{bial}
We now prove that the algebra and coalgebra structures on 
$A$ are compatible so that they give a bialgebra structure 
on $A$. 

Since 
\[
 \vep(\eta(\bfone)) = 
\vep( \bfone) = \vep(\base{0})= \bfone 
= \id (\bfone), \]
we have verified the commutativity of diagram~(\ref{ee}).

We also have 
\[ (\eta\otimes \eta) (\Delta_C(\bfone)) 
= (\eta\otimes \eta)(\bfone\otimes \bfone) 
= \base{0}\otimes \base{0}\]
and 
\[ \Delta (\eta (\bfone))=\Delta(\base{0}) 
= \sum_{k=0}^0 \sum_{i=0}^{0-k} (-\lambda)^k\base{i} 
    \otimes \base{0-k-i} 
=\base{0}\otimes \base{0}.\]
This proves the commutativity of diagram~(\ref{ed}).

We next prove the commutativity of diagram~(\ref{me}). 
We have 
\[ (\vep \otimes \vep)(\base{m}\otimes \base{n}) = 
\left \{ \begin{array}{ll} 
 \bfone \otimes \bfone, & (m,n)=(0,0), \\
\lambda\bfone\otimes \bfone, & (m,n)=(1,0),\\
\bfone \otimes \lambda \bfone, & (m,n)=(0,1), \\
\lambda\bfone \otimes \lambda\bfone, & (m,n)=(1,1), \\
0, & m\geq 2 {\rm\ or\ } n\geq 2. 
\end{array} \right . \]
So 
\begin{eqnarray*}
 \mu_C ((\vep\otimes \vep)(\base{m}\otimes \base{n}) )
&=& \left \{ \begin{array}{ll} 
 \bfone, & (m,n)=(0,0), \\
\lambda\bfone, & (m,n)=(1,0) {\rm\ or\ } (0,1),\\
\lambda^2 \bfone, & (m,n)=(1,1), \\
0, & m\geq 2 {\rm\ or\ } n\geq 2. 
\end{array} \right . 
\end{eqnarray*}
On the other hand, we have 
\begin{eqnarray*}
\vep(\mu(\base{m}\otimes \base{n})) 
&=& \vep \left ( \sum_{i=0}^m \lambda^i \binc{m+n-i}{m}\binc{m}{i} 
    \base{m+n-i} \right )\\
&=& \sum_{i=0}^m \lambda^i \binc{m+n-i}{m}\binc{m}{i} 
    \vep( \base{m+n-i} ). 
\end{eqnarray*}
So when $(m,n)=(0,0)$, we have 
\[ \vep(\mu(\base{0}\otimes \base{0})) 
= \vep (\base{0})= \bfone. \]
When $(m,n)=(0,1)$, we have 
\[
 \vep(\mu(\base{0}\otimes \base{1}))
= 
\sum_{i=0}^0 \lambda^i \binc{1-i}{0}\binc{0}{i} 
    \vep(\base{1-i}) 
= \vep(\base{1}) = \lambda\bfone. \]
By the commutativity of the multiplication $\mu$ in $A$, 
we also have 
\[ \vep(\mu(\base{1}\otimes \base{0}))= \lambda\bfone. \]
When $(m,n)=(1,1)$, we have 
\begin{eqnarray*}
 \vep(\mu(\base{1}\otimes \base{1})) &=& 
\sum_{i=0}^1 \lambda ^i \binc{2-i}{1}\binc{1}{i} 
    \vep(\base{2-i}) \\
&=& \lambda^0 \binc{2}{1}\binc{1}{0} \vep(\base{2}) 
+ \lambda \binc{1}{1}\binc{1}{1} \vep (\base{1}) \\
&=& \lambda^2 \bfone. 
\end{eqnarray*}

When $n\geq 2$, we have $m+n-i \geq 2$ for 
$0\leq i\leq m$. Thus 
$\vep( \base{m+n-i})=0$ and 
$\vep( \mu(\base{m}\otimes \base{n}))=0$. 
The same is true when $m\geq 2$ 
by the commutativity of $\mu$. 

Thus we have verified the commutativity of 
diagram~(\ref{me}).

The rest of this section is devoted to 
the verification the commutativity 
of diagram~(\ref{md}). 
In other words, we want to prove the identity
\margin{md1}
\begin{equation}
 \Delta (\mu(\base{m}\otimes \base{n})) 
=
    (\mu\otimes \mu)(\id\otimes \tau_{A,A}\otimes \id)
    (\Delta \otimes \Delta)(\base{m}\otimes \base{n}), 
\mlabel{md1}
\end{equation}
for any $m, n\geq 0$. 
The left hand side can be simplified as follows.  
\begin{eqnarray*}
\lefteqn{
\sum_{i=0}^m \sum_{k=0}^{m+n-i}\sum_{j=0}^{m+n-i-k}
    (-1)^k \lambda^{i+k} \bincc{m+n-i}{m}\bincc{m}{i}
    \base{j}\otimes \base{m+n-i-k-j} }\\
&=& 
\sum_{i=0}^{m+n} \sum_{k=0}^{m+n-i}\sum_{j=0}^{m+n-i-k}
    (-1)^k \lambda^{i+k} \bincc{m+n-i}{m}\bincc{m}{i}
    \base{j}\otimes \base{m+n-i-k-j} \\
&&  \left (\bincc{m}{i}=0 {\rm\ for\ } i>m \right ) \\
&=& 
\sum_{j=0}^{m+n} \sum_{k=0}^{m+n-j}\sum_{i=0}^{m+n-j-k}
    (-1)^k \lambda^{i+k} \bincc{m+n-i}{m}\bincc{m}{i}
    \base{j}\otimes \base{m+n-i-k-j} \\
 && {\rm( exchanging\ the\ first\ and\ third\ summations\ )}\\
&=& 
\sum_{j=0}^{m+n} \sum_{k=0}^{m+n-i}\sum_{\ell=0}^{m+n-j-k}
    (-1)^k \lambda^{m+n-j-\ell} \bincc{j+k+\ell}{m}
    \bincc{m}{m+n-j-k-\ell}
    \base{j}\otimes \base{\ell} \\
&&  ({\rm letting\ } i=m+n-j-k-\ell {\rm\ in\ the\ third\ sum} )\\
&=& 
\sum_{j=0}^{m+n} \sum_{\ell=0}^{m+n-j} \lambda^{m+n-j-\ell}  
\left ( \sum_{k=0}^{m+n-j-\ell} (-1)^k \bincc{j+k+\ell}{m}
    \bincc{m}{m+n-j-k-\ell} \right )
    \base{j}\otimes \base{\ell} \\
 && {\rm( exchanging\ the\ second\ and\ third\ summations\ )}.
\end{eqnarray*}

We next simplify the right hand side of equation~(\ref{md1}). 
Unwinding definitions we see that the right hand side is 

\begin{eqnarray*}
\lefteqn{\sum_{k=0}^m\sum_{i=0}^{m-k}\sum_{\ell=0}^n \sum_{j=0}^{n-\ell}
    \sum_{u=0}^i \sum_{v=0}^{m-k-i}
    (-1)^{k+\ell} \lambda^{k+\ell+u+v}}\\
 && \times \bincc{i+j-u}{i}\bincc{i}{u} \bincc{m+n-k-i-\ell-j-v}{m-k-i}
    \bincc{m-k-i}{v}
    \base{i+j-u}\otimes \base{m+n-k-i-\ell-j-v}.
\end{eqnarray*}
By substitutions 
\[ \left \{ \begin{array}{l}
b=i+j-u,\\
e=m+n-k-i-j-\ell-v
\end{array} \right . \]
where we treat $u$ and $v$ as the variables, we obtain 

\begin{eqnarray*}
\lefteqn{\sum_{k=0}^m\sum_{i=0}^{m-k}\sum_{\ell=0}^n \sum_{j=0}^{n-\ell}
    \sum_{b=j}^{i+j} \sum_{e=n-\ell-j}^{n+m-i-j-k-\ell}
    (-1)^{k+\ell} \lambda^{m+n-b-e}}\\
 && \times \bincc{b}{i}\bincc{i}{i+j-b} \bincc{e}{m-k-i}
    \bincc{m-k-i}{m+n-i-j-k-\ell-e}
    \base{b}\otimes \base{e}.
\end{eqnarray*}
Because of the nature of the summation limits, 
we cannot yet exchange the order of the summations as we did 
for the left hand side of equation~(\ref{md1}). But we 
have 

\begin{lemma}
\margin{md2}
\mlabel{md2}
\begin{eqnarray*}
\lefteqn{\sum_{k=0}^m\sum_{i=0}^{m-k}\sum_{\ell=0}^n \sum_{j=0}^{n-\ell}
    \sum_{b=j}^{i+j} \sum_{e=n-\ell-j}^{n+m-i-j-k-\ell}
    (-1)^{k+\ell} \lambda^{m+n-b-e}}\\
 && \times \bincc{b}{i}\bincc{i}{i+j-b} \bincc{e}{m-k-i}
    \bincc{m-k-i}{m+n-i-j-k-\ell-e}
    \base{b}\otimes \base{e}\\
\\
\lefteqn{=\sum_{k=0}^{m+n}\sum_{i=0}^{m+n}\sum_{\ell=0}^{m+n} 
\sum_{j=0}^{m+n}
    \sum_{b=0}^{m+n} \sum_{e=0}^{m+n}
    (-1)^{k+\ell} \lambda^{m+n-b-e}}\\
 && \times \bincc{b}{i}\bincc{i}{i+j-b} \bincc{e}{m-k-i}
    \bincc{m-k-i}{m+n-i-j-k-\ell-e}
    \base{b}\otimes \base{e}.
\end{eqnarray*}
\end{lemma}
\proofbegin
Note that we have $m,n,b,e\geq 0$ by assumption.  
Also for any integers $x,\ y$ with $x\geq 0$ and $y<0$ or 
$x\geq 0$ and $y>x$, we have 
$\bincc{x}{y}=0$. 
So 
\[ k>m  \Rightarrow m-k-i<0 \Rightarrow \bincc{e}{m-k-i}=0.\]
This shows that we can replace the first sum on the left hand 
side of the equation in the lemma by the first sum on the right hand side. 
Similarly, we have 
\[ \begin{array}{ccl}
i>m-k & \Rightarrow & m-k-i<0 \Rightarrow \bincc{e}{m-k-i}=0,\\
\ell>n &\Rightarrow & n-\ell-j<0 \Rightarrow 
\bincc{e}{n-\ell-j}=0,\\
j>n-\ell &\Rightarrow & n-\ell-j<0 \Rightarrow 
\bincc{e}{n-\ell-j}=0,\\
b<j &\Rightarrow & \bincc{b}{j}=0, \\
b>i+j &\Rightarrow & i+j-b<0 \Rightarrow \bincc{i}{i+j-b}=0,\\
e<n-j-\ell &\Rightarrow & \bincc{e}{n-j-\ell}=0, \\
e>m+n-i-j-k-\ell &\Rightarrow & 
    \bincc{m-k-i}{m+n-i-j-k-\ell-e}=0.
\end{array} \]
Considering in addition  
\[ \bincc{b}{i}\bincc{i}{i+j-b}=\bincc{b}{j}\bincc{j}{i+j-b}
\]
and 
\[ 
\bincc{e}{m-k-i}\bincc{m-k-i}{m+n-i-j-k-\ell-e}
=\bincc{e}{n-\ell-j}\bincc{n-\ell-j}{m+n-i-j-k-\ell-e},\]
we see that each of the other sums on the left hand side of the 
equation can be replaced by the corresponding sum on the right 
hand side of the equation. 
This proves the lemma. 
\proofend

Continuing with the proof of the commutativity of the 
diagram~(\ref{md}), we see that 
the limits of the sums on the right hand side of  the 
equation in Lemma~\ref{md2} 
are given by the same constants. 
Thus we can exchange the order of the summations and get 

\begin{eqnarray*}
\lefteqn{\sum_{k=0}^{m+n}\sum_{i=0}^{m+n}\sum_{\ell=0}^{m+n} 
\sum_{j=0}^{m+n}
    \sum_{b=0}^{m+n} \sum_{e=0}^{m+n}
    (-1)^{k+\ell} \lambda^{m+n-b-e}}\\
 && \times \bincc{b}{i}\bincc{i}{i+j-b} \bincc{e}{m-k-i}
    \bincc{m-k-i}{m+n-i-j-k-\ell-e}
    \base{b}\otimes \base{e}\\
\\
\lefteqn{= \sum_{b=0}^{m+n} \sum_{e=0}^{m+n}
\sum_{i=0}^{m+n}\sum_{j=0}^{m+n}\sum_{k=0}^{m+n} 
\sum_{\ell=0}^{m+n}
        (-1)^{k+\ell} \lambda^{m+n-b-e}}\\
 && \times \bincc{b}{i}\bincc{i}{i+j-b} \bincc{e}{m-k-i}
    \bincc{m-k-i}{m+n-i-j-k-\ell-e}
    \base{b}\otimes \base{e}.
\end{eqnarray*}
Now by the same argument as in the proof of Lemma~\ref{md2}, 
we obtain
\begin{eqnarray*}
\lefteqn{\sum_{b=0}^{m+n} \sum_{e=0}^{m+n}
\sum_{k=0}^{m+n}\sum_{i=0}^{m+n}\sum_{\ell=0}^{m+n} 
\sum_{j=0}^{m+n}
        (-1)^{k+\ell} \lambda^{m+n-b-e}}\\
 && \times \bincc{b}{i}\bincc{i}{i+j-b} \bincc{e}{m-k-i}
    \bincc{m-k-i}{m+n-i-j-k-\ell-e}
    \base{b}\otimes \base{e}\\
\\
\lefteqn{= \sum_{b=0}^{m+n} \sum_{e=0}^{m+n-b}
\left [\sum_{i=0}^{b}\sum_{j=b-i}^{b}\sum_{k=0}^{m-i} 
\sum_{\ell=0}^{m+n-i-j-k-e}
        (-1)^{k+\ell} \lambda^{m+n-b-e} \right .}  \\
 && \left. 
 \times \bincc{b}{i}\bincc{i}{i+j-b} \bincc{e}{m-k-i}
    \bincc{m-k-i}{m+n-i-j-k-\ell-e}\right ]
    \base{b}\otimes \base{e}.
\end{eqnarray*}
Comparing the right hand side of the above equation 
with the simplified form of the left hand side of 
equation~(\ref{md1}), 
we see that to prove equation~(\ref{md1}), we only need to prove 
 
\margin{md3}
\begin{equation}
\begin{array}{l}
\displaystyle{\sum_{k=0}^{m+n-b-e} (-1)^k \bincc{b+e+k}{m}
    \bincc{m}{m+n-b-e-k}}  \\
=
\displaystyle{\sum_{i=0}^{b}\sum_{j=b-i}^{b}\sum_{k=0}^{m-i} 
\sum_{\ell=0}^{m+n-i-j-k-e}
        \!\!\!\!\!\!\! (-1)^{k+\ell} 
\bincc{b}{i}\bincc{i}{i+j-b} \bincc{e}{m-k-i}
    \bincc{m-k-i}{m+n-i-j-k-\ell-e}}
\end{array}
\mlabel{md3}
\end{equation}
for all $m,n,b,e\geq 0$ with $b+e\leq m+n$. 

Using the substitutions 
$$ 
\left \{ \begin{array}{l}
j=b-i+c,\\
k=m-i-a,\\
\ell = n-b-c-e+i+a-d\end{array} \right . $$
on the right hand side of 
equation~(\ref{md3}) gives us 

\begin{eqnarray*}
\lefteqn{
\sum_{i=0}^{b}\sum_{j=b-i}^{b}\sum_{k=0}^{m-i} 
\sum_{\ell=0}^{m+n-i-j-k-e}
        (-1)^{k+\ell} 
\bincc{b}{i}\bincc{i}{i+j-b} \bincc{e}{m-k-i}
    \bincc{m-k-i}{m+n-i-j-k-\ell-e}}\\
&=&\sum_{i=0}^{b}\sum_{c=0}^{i}\sum_{a=0}^{m-i} 
\sum_{d=0}^{n-b-c-e+i+a}
        (-1)^{m+n-e-b-c-d} 
\bincc{b}{i}\bincc{i}{c} \bincc{e}{a}
    \bincc{a}{d}.
\end{eqnarray*}
Thus to prove equation~(\ref{md3}), and hence 
equation~(\ref{md1}), we only need to prove 
the following theorem.

\margin{ID1}
\mlabel{ID1}
\begin{theorem}
If $m, n, b$ and $e$ are nonnegative integers satisfying
\begin{displaymath}
m+n\geq b+e,
\end{displaymath}
then
\begin{displaymath}
\begin{array}{cc} 
\displaystyle{(-1)^{m+n-e-b} \sum_{k=0}^{m+n-b-e} (-1)^k
\left ( \begin{array}{cc} b+e+k\\ m \end{array}  \right )
\left ( \begin{array}{cc} m\\ m+n-b-e-k \end{array} \right )}
\ \ \ \ \ \ \\  \\ 
 \ \ \ \ \ \ 
\displaystyle{=\sum_{i=0}^{b} \sum_{c=0}^{i}
\sum_{a=0}^{m-i}\ \sum_{d=0}^{n-e-b+i-c+a}
(-1)^{c+d} \left ( \begin{array}{c} b\\i \end{array} \right )
\left ( \begin{array}{c} i\\ c \end{array} \right )
\left ( \begin{array}{c}  e\\ a\end{array}\right )
\left ( \begin{array}{c} a\\ d\end{array}\right ).}
\end{array}
\end{displaymath}
\end{theorem}

\proofbegin
By replacing $k$ by $m+n-b-e-k$
in the sum on the left hand side of the above equation 
and using the fact that~\cite{S-W}
\begin{displaymath}
\sum_{d=0}^j (-1)^d \left ( \begin{array}{cc} a\\d \end{array}
\right ) = (-1)^j \left ( \begin{array}{cc} a-1\\j \end{array}
\right )
\end{displaymath}
(Notice that this means that $\left ( \begin{array}{cc}-1\\j\end{array}
\right )=(-1)^j$),
the problem is reduced to showing that
\margin{id1}
\begin{equation}
\begin{array}{cc}
\displaystyle{\sum_{i=0}^{b}\sum_{c=0}^{i}\sum_{a=0}^{m-i} 
(-1)^{n-e-b+i+a} \bincc{b}{i}
\bincc{i}{c} \bincc{e}{a} \bincc{a-1}{n-e-b+i-c+a}} \\ \ \ 
\\
\ \ \ \ \ \ \ \ \displaystyle{=
\sum_{k=0}^{m+n-b-e}(-1)^k \bincc{m+n-k}{m}
\bincc{m}{k}}
\end{array}
\mlabel{id1}
\end{equation}
Using the classical summation of Vandermonde~\cite{S-W} 
\begin{displaymath}
\sum_{k=0}^m \left ( \begin{array}{cc}m\\ k \end{array}
\right ) \left ( \begin{array}{cc} n\\ i-k
\end{array}\right )=
\left ( \begin{array}{cc} m+n\\ i \end{array}
\right )
\end{displaymath}
in the summation on $c$ on the left hand side of (\ref{id1}), 
we find that (\ref{id1})
reduces to
\margin{id2}
\begin{equation}
\begin{array}{cc}
\displaystyle{\sum_{i=0}^b \sum_{a=0}^{m-i} (-1)^{n-e-b+i+a}
\left (\begin{array}{cc}b\\ i \end{array}\right )
\left ( \begin{array}{cc} e\\ a \end{array} \right )
\left ( \begin{array}{cc} a+i-1\\ n-e-b+i+a
\end{array}\right )}
\\ \ \ 
\\
\ \ \ \ \ \ \ \ \ \ \ \displaystyle{=
\sum_{k=0}^{m+n-b-e}(-1)^k \left (
\begin{array}{cc}m+n-k\\ m\end{array} \right )
\left ( \begin{array}{cc}m\\k \end{array}\right ).}
\end{array}
\mlabel{id2}
\end{equation}
Now set $T=a+i$. By the Vandermonde again, we find that
the left hand side of (\ref{id2}) becomes
\begin{eqnarray*}
\lefteqn{ \sum_{T\leq m} \sum_{i=0}^b (-1)^{n-e-b+T}\left ( \begin{array}{cc}
b\\ i\end{array}\right ) \left ( \begin{array}{cc}e\\T-i\end{array}
\right ) \left ( \begin{array}{cc} T-1\\ n-e-b+T \end{array}
\right )}\\
& =&  \sum_{T\leq m} (-1)^{n-e-b+T}\left ( \begin{array}{cc} b+e\\T
\end{array}\right ) \left ( \begin{array}{cc} T-1\\n-e-b+T
\end{array}\right ).
\end{eqnarray*}
Therefore, it suffices to prove, by letting $H=b+e$ in (\ref{id2}), the following
identity
\margin{id3}
\begin{equation}
\begin{array}{cc}
\displaystyle{\sum_{k=0}^{m+n-H}(-1)^k \left ( \begin{array}{cc}m+n-k\\m \end{array}
\right ) \left ( \begin{array}{cc}m\\ k
\end{array}\right )} \ \ \ \ \ \ 
\\ \ \ 
\\
\ \ \ \ \ \ \ \ \ \ \ \displaystyle{
=\sum_{T\leq m} (-1)^{n-H+T}\left ( \begin{array}{cc}H\\T
\end{array}\right ) \left ( \begin{array}{cc}T-1\\n-H+T
\end{array}
\right ).}
\end{array}
\mlabel{id3}
\end{equation}
For brevity, write (\ref{id3}) as
\margin{id4}
\begin{equation}
L(m,n)=R(m,n).
\mlabel{id4}
\end{equation}
Clearly we have the following 
\begin{eqnarray}
L(m,0)=L(0,n)=R(0,n)=R(m,0)=1, \ \ \ \ \ \ \ \ \ \ \ \ 
\\
R(m,n)-R(m-1,n)=(-1)^{n+m-H}\left ( \begin{array}{cc}H\\m\end{array}\right )
\left ( \begin{array}{cc}m-1\\n+m-H\end{array}\right ).
\end{eqnarray}
Therefore, we have 
\margin{id7}
\begin{equation}
\begin{array}{l}
\displaystyle{R(m,n)-R(m-1,n)-R(m,n-1)+R(m-1,n-1)}\\
\ \ \ \ \ \displaystyle{
=(-1)^{n+m-H}\left ( \begin{array}{cc}H\\m\end{array}\right )
\left ( \left ( \begin{array}{cc} m-1\\n+m-H\end{array}\right )+
\left ( \begin{array}{cc}m-1\\n-1+m-H\end{array}\right ) \right )}
\\
\\
\ \ \ \ \ 
\displaystyle{=(-1)^{n+m-H}\left (\begin{array}{cc}H\\m\end{array}\right )
\left ( \begin{array}{cc}m\\n+m-H \end{array}\right ).}
\end{array}
\mlabel{id7}
\end{equation}
Using the fact that
\begin{displaymath}
\left ( \begin{array}{cc} m+n-k\\m\end{array}\right )
\left ( \begin{array}{cc}m\\ k\end{array}\right )=
\left ( \begin{array}{cc}n\\ k\end{array}\right )
\left (\begin{array}{cc}m+n-k\\n\end{array}\right ),
\end{displaymath}
we find that
\begin{eqnarray*}
\lefteqn{L(m,n)-L(m-1,n)} \\
\\
&=& (-1)^{m+n-H}
\left ( \begin{array}{cc}H\\m \end{array}\right )\left (
\begin{array}{cc}m\\m+n-H\end{array} \right ) \\
\\
&& +\sum_{k=0}^{m+n-H-1}(-1)^k
\left (\begin{array}{cc}n\\k\end{array}\right )
\left ( \left ( \begin{array}{cc}m+n-k\\n \end{array}
\right )-\left (\begin{array}{cc} m+n-k-1\\n \end{array}
\right ) \right ) \\
\\
&=& (-1)^{m+n-H}
    \left ( \begin{array}{cc}H\\m \end{array} \right )
\left ( \begin{array}{cc}m\\m+n-H\end{array}\right )
\\
\\
 && +\sum_{k=0}^{m+n-H-1}(-1)^k
\left ( \begin{array}{cc}n\\ k\end{array}\right )
\left ( \begin{array}{cc}m+n-k-1\\m-k\end{array}\right ) \\
\\
&=& (-1)^{m+n-H} 
\left ( \begin{array}{cc}H\\m\end{array}\right )
\left (\begin{array}{cc}m\\m+n-H\end{array}\right )
\\
\\
 && +\sum_{k=0}^{m+n-H-1}(-1)^k\left ( \left (
\begin{array}{cc}n-1\\k\end{array}\right ) + \left ( \begin{array}{cc}
n-1\\k-1\end{array}\right ) \right ) \left ( \begin{array}{cc}
m+n-k-1\\m-k\end{array}\right )\\ 
\\
&=& (-1)^{m+n-H} 
    \left ( \begin{array}{cc}H\\m\end{array}\right )
\left ( \begin{array}{cc}m\\m+n-H\end{array}\right )
+L(m,n-1)
\\
\\
&& +\sum_{k=0}^{m+n-H-2}(-1)^{k+1}\left (
\begin{array}{cc}n-1\\k\end{array}\right )
\left ( \begin{array}{cc}m+n-k-2\\m-1-k\end{array}\right )\\
\\
&=& (-1)^{m+n-H} 
    \left ( \begin{array}{cc}H\\m\end{array}\right )
\left ( \begin{array}{cc}m\\m+n-H\end{array}\right )
+L(m,n-1)-L(m-1,n-1).
\end{eqnarray*}
Therefore, 
\margin{id8}
\begin{equation}
\begin{array}{l}
\displaystyle{L(m,n)-L(m-1,n)-L(m,n-1)+L(m-1,n-1)}\\ 
\\
\ \ \ \ \ \
 \displaystyle{=(-1)^{n+m-H}\left 
(\begin{array}{cc}H\\m\end{array}\right ) \left ( 
\begin{array}{cc}m\\n+m-H \end{array}\right ).} 
\end{array}
\mlabel{id8}
\end{equation}
Thus we see that (\ref{id7}) and (\ref{id8}) show that $L(m,n)$ and
$R(m,n)$ satisfy the same bilinear recurrence.
Since they have the same initial values, we have that
$L(m,n)=R(m,n)$ for all nonnegative $n$ and $m$.
\proofend

\subsection{Existence of an Antipode}
\mlabel{hopf}
We now show that the linear map $S$ defined in  
Section~\ref{intro} is an antipode on the 
bialgebra $\cala_\lambda$, thus making the bialgebra 
into a Hopf algebra. 

Since $A$ is commutative, we only need to prove 
\[ \mu\circ (S\otimes \id)\circ \Delta = \eta\circ \vep.\]
From the definitions of $\vep$ and $\eta$, we have 
\[ \eta\circ \vep (\base{n}) = \left \{ 
\begin{array}{ll} 
\base{0}, & n=0,\\
\lambda\base{0}, & n=1,\\
0, & n>1.
\end{array}
\right . 
\]
Thus to prove that $S$ is an antipode on $A$, we only need to 
show  
\margin{id10}
\begin{equation}
\ \mu(S\otimes \id(\Delta(\base{n}))) = 
    \left \{ 
    \begin{array}{ll} 
    \base{0}, & n=0,\\
    \lambda\base{0}, & n=1,\\
    0, & n>1.
    \end{array}
    \right . 
\mlabel{id10}
\end{equation}
Recall that we have defined the $C$-linear map 
$S:A\to A$ by 
\[
S(\base{n})= 
  (-1)^n \sum_{v=0}^n \bincc{n-3}{v-3} \lambda^{n-v} \base{v}.
\]
Using this and the definitions of $\mu$ and $\Delta$, we have 
$$
\mu(S\otimes \id(\Delta(\base{n}))) 
= \sum_{k=0}^n\sum_{i=0}^{n-k} \sum_{v=0}^i \sum_{\ell=0}^v
    (-1)^{k+i}  \lambda^{k+i-v+\ell}   \bincc{i-3}{v-3} 
     \bincc{n-k-i+v-\ell}{v}
    \bincc{v}{\ell} \base{n-k-i+v-\ell}.$$
With a change of variable
$ \ell =n-k-i+v-w\ (w=n-k-i+v-\ell))$, we get  
$$\mu(S\otimes \id(\Delta(\base{n}))) 
= \sum_{k=0}^n\sum_{i=0}^{n-k} \sum_{v=0}^i 
    \sum_{w=n-k-i}^{n-k-i+v}
    (-1)^{k+i}  \lambda^{n-w}   \bincc{i-3}{v-3} 
     \bincc{w}{v}
    \bincc{v}{n-k-i+v-w} \base{w}.$$

As in the proof of diagram~(\ref{md}), we want  
constant limits for the summations. 
Since $\bincc{x}{y}=0$ for integers $x,\ y$ with  
$x\geq 0$ and either $y<0$ or $y>x$, 
we have 

\[ \begin{array}{ll}
&i>n-k  \Rightarrow   n-k-i<0 \\
\\
\Rightarrow &
    \left \{ \begin{array}{ll} 
    n-k-i+v-w<0 \Rightarrow \bincc{v}{n-k-i+v-w}=0, 
     & {\rm\ if\ } v\leq w, \\
    \bincc{w}{v} =0, &
     {\rm\ if\ } v>w. 
    \end{array}   \right .
\end{array}
\]

\[\begin{array}{ccl}
w< n-k-i &  \Rightarrow &  n-k-i-w>0 \\
\\
    \Rightarrow  n-k-i+v-w>v 
    & \Rightarrow & \bincc{v}{n-k-i+v-w} =0, 
\end{array}\]
\[\begin{array}{ccl}
w>n-k-i+v & \Rightarrow & n-k-i+v-w<0  \\
\\    \Rightarrow    \bincc{v}{n-k-i+v-w}=0. && 
\end{array}
\]
Therefore the last nested sum can be replaced by 
\begin{eqnarray*}
&&\sum_{k=0}^n\sum_{i=0}^{n} \sum_{v=0}^i 
    \sum_{w=0}^{n}
    (-1)^{k+i}  \lambda^{n-w}   \bincc{i-3}{v-3} 
     \bincc{w}{v}
    \bincc{v}{n-k-i+v-w} \base{w}\\
&=& \sum_{w=0}^n \lambda^{n-w} 
    \left [ \sum_{k=0}^{n} \sum_{i=0}^n 
    \sum_{v=0}^{i}
    (-1)^{k+i}     \bincc{i-3}{v-3} 
     \bincc{w}{v}
    \bincc{v}{n-k-i+v-w} \right ] \base{w}.
\end{eqnarray*}
Thus to prove equation~(\ref{id10}) we only need to 
prove 
\margin{ID2}
\begin{theorem}
\mlabel{ID2}
For any integers $n\geq w\geq 0$, we have 
\margin{id11}
\begin{equation}
\begin{array}{l}
\displaystyle{ \sum_{k=0}^{n} \sum_{i=0}^n 
    \sum_{v=0}^{i}
    (-1)^{k+i}  \bincc{i-3}{v-3} 
     \bincc{w}{v}
    \bincc{v}{n-k-i+v-w}  }\\
\\
\ \ \ = \left \{ \begin{array}{ll} 
    1, & {\rm\ if\ } (n,w)=(0,0) {\rm\ or\ } (1,0), \\
    0, & {\rm\ otherwise\ }. 
    \end{array} \right .
\end{array}
\mlabel{id11}
\end{equation}
\end{theorem}

\proofbegin
The proof of this identity requires only three facts~\cite{Rio}:
\margin{id12-14}
\begin{eqnarray}
\left (\begin{array}{cc}-A\\i\end{array}\right ) = (-1)^i \left (
\begin{array}{cc}A+i-1\\i \end{array}\right ),
\mlabel{id12}\\
\sum_{j\geq 0} \left ( \begin{array}{cc}A\\j\end{array}\right )
\left (\begin{array}{cc}B\\C-j\end{array}\right )=
\left (\begin{array}{cc}A+B\\C\end{array}\right ),
\mlabel{id13}\\
\sum_{j=0}^{N} (-1)^j \left ( \begin{array}{cc}N\\j \end{array}
\right )= \left \{ \begin{array}{ll}
  1, & {\rm if\ } N=0,\\
   0, & {\rm otherwise.\ } 
\mlabel{id14}
\end{array}     \right .
\end{eqnarray}
Hence we have that
\begin{eqnarray*}
\lefteqn{\sum_{k=0}^{n}\sum_{i=0}^{n}\sum_{v=0}^{i}(-1)^{k+i}
\left (\begin{array}{cc}i-3\\v-3\end{array}\right )
\left ( \begin{array}{cc}w\\v\end{array}\right )
\left ( \begin{array}{cc}v\\n-k-i+v-w\end{array}\right )}\\
& =& \sum_{k=0}^{n}
\sum_{v=0}^{n}\sum_{i=v}^{n}
(-1)^{k+i}
\left (\begin{array}{cc}i-3\\v-3\end{array}\right )
\left ( \begin{array}{cc}w\\v\end{array}\right )
\left (\begin{array}{cc}v\\n-k-i+v-w\end{array}\right )\\
&=&
\sum_{k=0}^{n}\sum_{v=0}^{n}\sum_{i=0}^{n-v}(-1)^{k+i+v}
\left ( \begin{array}{cc}i+v-3\\v-3\end{array}\right )
\left ( \begin{array}{cc}w\\v\end{array}\right )
\left ( \begin{array}{cc}v\\n-k-i-w\end{array}\right ).
\end{eqnarray*}
Now by (\ref{id12}) and (\ref{id13}) we find that this sum equals
\begin{eqnarray*}
\lefteqn{\sum_{k=0}^{n} \sum_{v=0}^{n} (-1)^{k+v}\left (
\begin{array}{cc}w\\v \end{array}\right )
\sum_{i=0}^{n-v}\left ( \begin{array}{cc}2-v\\i\end{array}\right )
\left (\begin{array}{cc} v\\n-k-i-w\end{array}\right )}\\
\ \ \ \ \ \ \ \ &=&
\sum_{k=0}^{n}\sum_{v=0}^{n} (-1)^{k+v}\left (
\begin{array}{cc}2\\n-k-w\end{array}\right )
\left ( \begin{array}{cc}w\\v\end{array}\right ).
\end{eqnarray*}
By (\ref{id14}) we find that this sum equals
\begin{displaymath}
\begin{array}{ll}
\ \ \ \left \{ \begin{array}{ll} 0, & {\rm if \ } w>0,\\
   \displaystyle{\sum_{k=0}^{n}(-1)^k 
    \left ( \begin{array}{cc}2\\n-k\end{array}\right )},
& {\rm if\ } w=0. 
\end{array} \right .\\
\ \ \ \\
\ \ \ = \left \{ \begin{array}{ll} 0, & {\rm if \ } w>0,\\
                                   1, & {\rm if \ } w=0 {\rm \ and\ } n=0,\\
                                   1, & {\rm if \ } w=0 {\rm \ and \ } n=1,\\
                                   0, & {\rm if \ } n>1. 
\end{array} \right . 
\end{array}
\end{displaymath}

\proofend

\medskip
This completes the proof of Theorem 1.1.

\section{Isomorphisms between $\lambda$-divided power Hopf algebras}
\mlabel{sec:iso}

We will prove the three statements in Theorem~\ref{thm:iso} by proving 
Proposition~\ref{pp:iso1}, \ref{pp:iso2} and \ref{pp:iso3}.  
To distinguish elements in $\lambda$-divided power Hopf algebras 
with various values of $\lambda$, we write 
$$A_\lambda=\bigoplus_{n\geq 0} Ca_{\lambda,n}$$
where $\{a_{\lambda,n}\}$ is the standard basis. 

\begin{prop}
Let $\lambda,\, \nu$ be in $C$. 
If $(\lambda)=(\nu)$, then $\cala_\lambda$ and $\cala_\nu$ are 
isomorphic Hopf algebras. 
\mlabel{pp:iso1}
\end{prop}
\begin{proof}
Suppose $(\lambda)=(\nu)$. Then $\lambda=\omega\nu$ for a unit 
$\omega$ in $C$. 
Define 
$$\varphi: A_\lambda \to A_\nu,\ 
a_{\lambda,n} \mapsto \omega^na_{\nu,n}.$$
Then it is straightforward to verify that $\varphi$ is a 
homomorphism of Hopf algebras. More precisely, we have 
\begin{eqnarray*}
\varphi \circ \mu_\lambda &=&\mu_\nu \circ (\varphi\otimes\varphi),\\
(\varphi\otimes \varphi) \circ \Delta_\lambda &=&\Delta_\nu \circ 
\varphi,\\
\varphi \circ S_\lambda &=&S_\nu \circ \varphi.
\end{eqnarray*}
It is also easily seen that $\varphi$ has inverse given by 
sending $a_{\nu,n}$ to $\omega^{-n} a_{\lambda,n}$. 
\end{proof}   

\begin{prop}
Let $C$ be a $\QQ$-algebra, and let $\lambda,\, \nu$ be in $C$.  
Then $\cala_\lambda$ and $\cala_\nu$ are 
isomorphic Hopf algebras if and only if $(\lambda)=(\nu)$.
\mlabel{pp:iso2}
\end{prop}
\begin{proof}
By Proposition~\ref{pp:iso1}, we only need to prove that 
if $\cala_\lambda$ and $\cala_\nu$ are 
isomorphic Hopf algebras, then $(\lambda)=(\nu)$.

%We will prove by contradiction. 
Suppose $C$ is a $\QQ$-algebra  
%$(\lambda)\neq (\nu)$ 
and $\cala_\lambda\cong 
\cala_\nu$. 
Let $\varphi: \cala_\lambda\to \cala_\nu$ be a Hopf algebra 
isomorphism. By \cite[Proposition 3.2]{Gu}, both $A_\lambda$ and 
$A_\nu$ are isomorphic to $C[x]$ as $C$-algebras with  
$a_{\lambda,1}$ and $a_{\nu,1}$ as the generators. It follows 
that $\varphi(a_{\lambda,1})=\sum_{i\geq 1} c_i a_{\nu,i}$ with 
$c_1\in C^*$. Likewise, 
$\varphi^{-1}(a_{\nu,1})=\sum_{i\geq 1} d_i a_{\lambda,i}$ with 
$d_1\in C^*$. 
Then we have 
\begin{eqnarray*}
((\varphi\otimes \varphi)\circ \Delta_\lambda)(a_{\lambda,1}) 
&=&(\varphi\otimes \varphi)(1\otimes 
a_{\lambda,1}+a_{\lambda,1}\otimes 1 -\lambda \otimes 1)\\
&=& 1\otimes \left (\sum_{i\geq 1} c_i a_{\nu,i}\right)+\left (\sum_{i\geq 
1}c_i a_{\nu,i}\right )\otimes 1 -\lambda\otimes 1\\
&=& -\lambda\otimes 1 + {\rm\ higher\ degree\ terms}
\end{eqnarray*}
and 
\begin{eqnarray*}
(\Delta_\nu\circ \varphi)(a_{\lambda,1})&=& 
\Delta_\nu\left (\sum_{i\geq 1} c_i a_{\nu,i}\right )\\
&=& c_1(1\otimes a_{\nu,1}+a_{\nu,1}\otimes 1-\nu \otimes 
1) \\ &+&c_2(1\otimes a_{\nu,2}+a_{\nu,1}\otimes 
a_{\nu,1}+a_{\nu,2}\otimes 1-\nu(1\otimes 
a_{\nu,1}+a_{\nu,1}\otimes 1)+\nu^2\otimes 1)\\
&+& 
\ldots\\
&=& -\nu\left (\sum_{i\geq 1} (-\nu)^{i-1}c_i\right )(1\otimes1)+ {\rm\ 
higher\ degree\ terms}.
\end{eqnarray*}
Here the degree is the natural one given by the standard basis 
$a_{\nu,i}$. That is, the degree of $a_{\nu,i}\otimes a_{\nu,j}$ 
is $(i,j)$ with the lexicographic order.
Since $\varphi$ is a Hopf algebra isomorphism, we have 
$$ \lambda=\nu\left (\sum_{i\geq 1} (-\nu)^{i-1}c_i\right ).$$
So $C\lambda \subseteq C\nu$. 
Considering $\varphi^{-1}$, we similarly get 
$$\nu=\lambda\left (\sum_{i\geq 1} (-\lambda)^{i-1} d_i\right ).$$
So $C\nu \subseteq C\lambda$. 
Then $C\nu=C\lambda$. %This is a contradiction. 
\dfootnote{Note that I only ask that coefficient of 1 in 
$\varphi(a_{\lambda,1})$ and $\varphi^{-1}(a_{\nu,1})$ be zero. 
Can this be used for a more general case?}
\end{proof}

For future reference, we record the following easy corollary of 
Proposition~\ref{pp:iso2}. 

\begin{coro}
Let $C$ be a $\QQ$-algebra and let $\nu$ be a unit in $C$. 
Then $\cala_\nu\not\cong \cala_0$. 
\mlabel{co:q}
\end{coro}

\begin{prop}
Let $\lambda$ and $\nu$ be in $C$. If either 
\begin{enumerate} 
\item $\nu$ is not contained in 
$\sqrt{(\lambda)}$, or 
\item 
there is a prime number $p$ such that 
$\nu^{p-1}$ is not contained in $(\lambda,p)$, 
\end{enumerate}
then $\cala_\lambda$ 
is not isomorphic to $\cala_\nu$ as $C$-algebras
\mlabel{pp:iso3}
\end{prop}
The third statement of Theorem~\ref{thm:iso} follows from the 
first case of the proposition. To display the utility of the second 
case, we provide the following examples. 

\begin{exam}
\begin{enumerate}
\item
Let $C=\FF_p[x]/(x^r)$ with $r\geq p$. Then 
$\cala_{\bar{x}} \not\cong \cala_0$. Here $\bar{x}$ is the image of $x$ in 
$\FF_p[x]/(x^r)$.
\item
Let $C=\ZZ[x]/(x^r)$ with $r\geq 2$. Then 
$\cala_{\bar{x}} \not\cong \cala_0$ $($taking $p=2$$)$. 
\end{enumerate}
\mlabel{ex:iso}
\end{exam}

\noindent
{\em Proof of Proposition~\ref{pp:iso3}.}
To consider the case when $\nu$ is not contained in $\sqrt{(\lambda)}$, 
we begin with a special situation. 

\begin{lemma}
Let $C$ be reduced with characteristic $p$ for a prime number 
$p$. Let $\nu$ be a unit in $C$. Then $\cala_0\not\cong \cala_\nu$. 
\mlabel{lem:01}
\end{lemma}
\begin{proof}
By \cite[Theorem 4.8]{Gu0}, the nilradical $N(A_0)$ of $A_0$ is 
$\bigoplus_{n\geq 1} C a_{0,n}.$
By \cite[Lemma 4.9]{Gu0}, we have $N(A_0)^p=0$. 
So we only need to show that $A_\nu$ has no non-zero nilpotent 
element $x$ with $x^p=0$. 

Let there be such an element $x\in A_\nu$. Write 
$x=\sum_{i=k}^\infty c_i a_{\nu,i}$ with $c_k\neq 0$. 
Then by the product formula in $A_\nu$ (see (ii) in the introduction), 
$$x^p=c_k^pa_{\nu,k}^p+{\rm\ a\ term\ in} \bigoplus_{i=k+1}^\infty Ca_{\nu,i}$$
and 
$$a_{\nu,k}^p=\nu^{\,(p-1)k} a_{\nu,k} +{\rm\ a\ term\ in} \bigoplus_{i=k+1}^\infty Ca_{\nu,i}.$$
Therefore,
$$x^p=c_k^p\nu^{p-1} a_{\nu,k}+{\rm\ a\ term\ in} \bigoplus_{i=k+1}^\infty Ca_{\nu,i}.$$
So we must have $c_k^p\nu^{p-1}=0$. Since $\nu$ is a unit, we have $c_k^p=0$. 
Since $C$ is reduced, we have $c_k=0$. This is a contradiction. 
\end{proof}

Continuing with the proof of Proposition~\ref{pp:iso3}, we now suppose  
$\nu\not\in \sqrt{(\lambda)}$, and $\cala_\lambda \cong \cala_\nu$.
Then $\nu$ is not nilpotent and $(\lambda)\cap \{\nu^n\}$ is empty. 
By \cite[Example 4]{Ma}, the image $\lambda/1$ of $\lambda$ in 
the localization $\{\nu^n\}^{-1}C$ is not a unit. 
Let $P$ be a maximal ideal of $\{\nu^n\}^{-1}C$ containing $\lambda/1$. 
Then the image $\widetilde{\lambda}$ of $\lambda/1$ (resp. of $\widetilde{\nu}$ of 
$\nu/1$) in the field $\{\nu^n\}^{-1}C/P$ is 
zero (resp. a unit). From $\cala_\lambda \cong \cala_\nu$, we have 
$$(\{\nu^n\}^{-1}C/P)\otimes_C A_\lambda\cong 
(\{\nu^n\}^{-1}C/P)\otimes_C A_\nu.$$
That is 
$\cala_{\widetilde{\lambda}}\cong 
\cala_{\widetilde{\nu}}.$
This is a contradiction by Corollary~\ref{co:q} (when the characteristic of the field 
$\{\nu^n\}^{-1}C/P$ is 0) and 
Lemma~\ref{lem:01} (when the characteristic is a prime number $p$).

\medskip
We next consider the case when there is a prime number $p$ such that $\nu^{p-1}$ is not 
contained in the ideal $(\lambda,p)$ of $C$. We again start with a 
special situation. 

\begin{lemma}
Let $C$ be of characteristic $p$. If $\nu^{p-1}$ is not zero, then 
$\cala_\nu$ is not isomorphic to $\cala_0$. 
\mlabel{lem:p1}
\end{lemma}
\begin{proof}
Suppose $\cala_\nu$ is isomorphic to $\cala_0$. Let $\varphi:A_0 \to 
A_\nu$ be an isomorphism of $C$-algebras. By \cite[Theorem 
4.8]{Gu0}, the nilradical $N(A_0)$ of $A_0$ is 
$N(C)\bigoplus \left (\bigoplus_{n\geq 1} C a_{0,n}\right ).$
On the other hand, by the product formula for $a_{\nu,n}\in A_\nu$ and the fact 
that $\bigoplus_{n\geq 1} Ca_{\nu,n}$ is an ideal of $A_\nu$, we have 
$N(A_\nu)=N(C)\bigoplus L$ for an ideal $L$ of $A_\nu$ contained 
in $\bigoplus_{n\geq 1} C a_{\nu,n}.$
Then $\varphi$ induces an isomorphism between the 
$C/N(C)$-algebras $A_0/N(A_0)\cong C/N(C)$ and $A_\nu/N(A_\nu) 
\cong C/N(C) \bigoplus \left ( \oplus_{n\geq 1} Ca_{\nu,n} 
\right)/L$. It follows that 
$L =\bigoplus_{n\geq 1} Ca_{\nu,n}.$
Thus $a_{\nu,1}$ is nilpotent. By \cite[Lemma 4.9]{Gu0}, 
$N(A_0)^p$ is contained in $N(C)^p$.
So we must have $a_{\nu,1}^p\in N(C)^p$. In particular, 
$a_{\nu,1}^p\in C$. 
But this cannot be true since 
$$a_{\nu,1}^p=\nu^{p-1} a_{\nu,1} +{\rm\ a\ term\ in} \bigoplus_{i=k+1}^\infty Ca_{\nu,i}$$
and $\nu^{p-1}\neq 0$.
\end{proof}

In the general case, consider $\overline{C}=C/(\lambda,p)$ and let 
$\bar{\lambda}$ (resp. $\bar{\nu}$) be the image of $\lambda$ (resp. 
$\nu$) in $\overline{C}$. Then 
by Lemma~\ref{lem:p1}, $\cala_{\bar{\nu}}\not\cong \cala_{\bar{\lambda}}$. 
So $\cala_\nu \not\cong \cala_\lambda$.
$\square$

\addcontentsline{toc}{section}{\numberline {}References}


\begin{thebibliography}{abcdsfgh}
\bibitem{Ba} G. Baxter, {An analytic problem whose solution
    follows from a simple algebraic identity,}
    {\em Pacific J. Math.} {\bf 10} (1960), 731-742.
     
\bibitem{B-O} P. Berthelot and A. Ogus, 
    {Notes on Crystalline Cohomology,}
    Princeton University Press, 1978. 

\bibitem{Ca} P. Cartier, {On the structure of free Baxter algebras,}
    {\em Adv. in Math.} {\bf 9} (1972), 253-265.
    
\bibitem{Ch} K.T. Chen, Integration of paths,
geometric invariants and a generalized Baker-Hausdorff
formula, {\em Ann. of Math.} {\bf 65} (1957), 163-178.

\bibitem{Gu0} L. Guo, {Properties of Baxter algebras,} 
    {\em Adv. in Math.} {\bf 151} (2000), 346-374. 
    
\bibitem{Gu} L. Guo, 
    {Baxter algebras and the umbral calculus,} 
    {\em Adv. in Appl. Math.} {\bf 27} (2001), 405-426. 

\bibitem{Gu1} L. Guo, {Baxter algebras and differential algebras,}
    In: Differential Algebra and Related Topics, World Scientific 
    (2002), 281-305. 

\bibitem{G-K1} L. Guo and W. Keigher, {Baxter algebras and
    shuffle products,}
    {\em Adv. in Math.} {\bf 150} (2000), 117-149.

\bibitem{G-K2} L. Guo and W. Keigher, {On free Baxter algebras:
    completions and the internal construction,}
    {\em Adv. in Math.} {\bf 151} (2000), 101-127.   

\bibitem{Ke} W. Keigher, { On the ring of Hurwitz series,}
    {\em Comm. Algebra} {\bf 25} (1997), 1845-1859.

\bibitem{Ma} H. Matsumura, {Commutative Ring Theory,}
    Cambridge University Press, 1994. 

\bibitem{N-S} W. Nichols and M. Sweedler, 
    {Hopf algebras and combinatorics, } In: 
    Umbral Calculus and Hopf Algebras, 
    Contemporary Mathematics {\bf 6}, 
    Amer. Math. Soc. (1982), 49-84. 

\bibitem{Re} R. Ree, { Lie elements and an algebra
    associated with shuffles,}
    {\em Ann. Math.} {\bf 68} (1958),
    210-220. 

\bibitem{Rio} J. Riodan, {Combinatorial Identities,} 
    Wiley, 1968. 

\bibitem{Rom} S. Roman, { The Umbral Calculus, }
    Academic Press, Orlando, FL, 1984. 

\bibitem{R-R} S. Roman and G.-C. Rota,  
    {The umbral calculus,} {\em Adv. Math.} {\bf 27}(1978), 
    95--188. 

\bibitem{Ro1} G.-C. Rota, {Baxter algebras and combinatorial
    identities I, II,} {\em Bull. Amer. Math. Soc.} {\bf 75} 
    (1969), 325--329, 330--334.
    
\bibitem{Ro2} G.-C. Rota, {Baxter operators, an introduction,}
    In: Gian-Carlo Rota on Combinatorics, Introductory Papers
    and Commentaries, Joseph P.S. Kung, Editor,
    Birkh\"{a}user, Boston, 1995.

\bibitem{Ro3} G.-C. Rota, {Ten mathematics problems I will never
    solve,} Invited address at the joint meeting of the
    American Mathematical Society and the Mexican Mathematical
    Society, Oaxaca, Mexico, December 6, 1997.
    DMV Mittellungen Heft 2, 1998, 45--52. 
        

\bibitem{S-W} 
D. Stanton and D. White, Constructive Combinatorics,
Springer-Verlag, New York, 1986. 

\end{thebibliography}
\end{document}